\documentclass[twoside]{svjour2}
\usepackage{tabularx}
\usepackage{float}
\usepackage{amsmath}
\usepackage{amsfonts}
\usepackage{amssymb}
\usepackage{amstext}
\usepackage{amscd}
\usepackage{graphicx}
\usepackage{enumerate}
\usepackage{pxfonts}

\def\EMn#1#2{#1^{({\rm EM}),#2}}
\def\CB3#1#2{#1^{({\rm cub3}),#2}}

\def\OD2#1#2{#1^{({\rm ord}\,2),#2}}
\def\ODp#1#2{#1^{({\rm ord}\,p),#2}}

\def\NCUB#1#2{#1^{({\rm New}),#2}}
\def\MC{{\rm MC}}

\def\VAR{{\rm Var}}

\begin{document}
  \title{Weak approximation of stochastic differential equations
    and application to derivative pricing
    \thanks{
      This research was partially supported by
      the Japanese Ministry of Education, Science, Sports and Culture,
      Grant-in-Aid for Scientific Research (C),
      15540110, 2003.
    }
  }
\author{Syoiti Ninomiya \and
  Nicolas Victoir}
\authorrunning{Ninomiya and Victoir}
\institute{Syoiti Ninomiya \at
  Center for Research in Advanced Financial
  Technology, Tokyo Institute of Technology, 2-12-1 Ookayama, Meguro-ku,
  Tokyo 152-8552 Japan \\
  \email{ninomiya@craft.titech.ac.jp}
  \and
  Nicolas Victoir \at
  Mathematical Institute, 24-29 St~Giles', Oxford, OX1 3LB, UK \\
  \email{victoir@maths.ox.ac.uk}
  \emph{Present email address:} \email{victoir@gmail.com}
}
\date{}


\maketitle

\begin{abstract}
  The authors present a new simple algorithm to approximate weakly stochastic
  differential equations
  in the spirit of~\mbox{\cite{kusuoka:2001aprx}\cite{LyonsVictoir:2002}.}
  They apply it to the
  problem of pricing Asian options under the Heston stochastic volatility
  model, and compare it with other known methods.
  It is shown that the
  combination of the suggested algorithm and quasi-Monte Carlo methods
  makes computations extremely fast.\\
  \\
\emph{2000 Mathematics Subject Classification.} 65C30, 65C05.\\
\keywords{
 Heston model,
 numerical methods for stochastic differential equations,
 mathematical finance, quasi-Monte Carlo method}
\end{abstract}
\section{Introduction}
\subsection{The Problem and its Motivation}\label{ss:problem}
We consider a stochastic differential equation written in the Stratonovich
form
\begin{equation}
  \label{StSDE}
  \begin{split}
    Y(t,x)&=x+\int_{0}^{t}V_{0}\left( Y(s,x)\right)\,ds
    +\sum_{i=1}^{d}\int_{0}^{t}V_{i}\left(Y(s,x)\right)\circ dB_{s}^{i},\\
    V_j &\in C_b^\infty\left(\mathbb{R}^N;\mathbb{R}^N\right),
  \end{split}
\end{equation}
where $B=\left(B^{1},\cdots ,B^{d}\right) $ is a standard
Brownian motion, and $C_b^\infty\left(\mathbb{R}^N;\mathbb{R}^N\right)$
denotes the set of $\mathbb{R}^N$-valued smooth functions defined over
$\mathbb{R}^N$ whose derivatives of any order are bounded.
In particular, we will use the classical notation
$Vf(x)=\sum_{i=1}^{N}V^{i}\left( x\right)
\left({\partial f}/{\partial x_{i}}\right)
\left( x\right)$
for $V\in C_{b}^{\infty }(\mathbb{R}^{N};\mathbb{R}^{N})$
and $f$ a differentiable function from $\mathbb{R}^{n}$ into $\mathbb{R}.$
This stochastic differential equation can be written in
It\^{o} form:
\begin{equation*}
Y(t,x)=x+\int_{0}^{t}\tilde{V}_{0}\left(Y(s,x)\right)\,ds
+\sum_{i=1}^{d}\int_{0}^{t}V_{i}\left(Y(s,x)\right)\,dB_{s}^{i},
\end{equation*}
where 
\begin{equation*}
\tilde{V}_{0}^{i}\left( y\right) =V_{0}^{i}\left(y\right) +\frac{1}{2}
\sum_{j=1}^dV_{j}V_{j}^{i}\left( y\right).
\end{equation*}%

Now, given a
function
$f$ with some regularity,
how can one approximate efficiently
$E\left[ f\left(Y(1,x)\right) \right] $?
It is equivalent to the following deterministic
problem:
if $L$ is the differential operator
$V_{0}+(1/2)\sum_{i=1}^{d}V_{i}^{2}$ and $u$
is the solution of the heat equation
\begin{equation*}
\frac{\partial u}{\partial t}\left(t, x\right) = Lu,\quad
u\left(0, x\right) = f(x),
\end{equation*}
how does one approximate $u\left( 1,x\right)$ (which is equal to
$E\left[f\left(Y(1,x)\right) \right]$
by Feynman-Kac theorem~\cite{ikeda-watanabe}).

This problem has had a lot of attention because of its practical
importance:\ it gives the evolution of the temperature in some media, and
also represents price of financial derivatives under stochastic financial
models such as Black-Scholes~\cite{BlackScholes:1973}.
\par
Non-probabilistic methods to solve the PDE (such as finite difference
methods) seem to only work well when $L$ is elliptic and in low dimension.
We refer to~\cite{LapeyrePardouxSentis}
for
a more detailed discussion on the subject.
We will focus in this
paper on probabilistic methods.
\subsection{Notation}
If $V$ is a smooth vector field, i.e. an element of $C_{b}^{\infty }\left( 
\mathbb{R}^{N};\mathbb{R}^{N}\right) $, $\exp \left( V\right)x$ denotes
the solution at time $1$ of the ordinary differential equation 
\begin{equation*}
\frac{dz_{t}}{dt}=V\left( z_{t}\right) ,\quad z_{0}=x.
\end{equation*}
For $x\in\mathbb{R}$, $\lfloor x\rfloor$
denotes the integer part of $x$.
For a random variable $X$, $\VAR[X]$ denotes the variance of $X$.

\subsection{Probabilistic Methods}
\subsubsection{Order 1}

The most popular probabilistic method to approximate $E\left[ f\left(
Y(1,x)\right) \right] $
is called the Euler-Maruyama method~\cite{KloedenPlaten:1999}.
We first fix $n$ independent $d$-dimensional random
variables $Z_{1},\cdots ,Z_{n}$
such that, if $X$ denotes
a standard normal random variables, 
\begin{equation}\label{eq:cub3cond}
  E\left[ p\left( Z_{k}\right) \right] =E\left[ p\left( X\right) \right]
\end{equation}
for all polynomial $p$ of degree less than or equal to $3$.
Then
one defines recursively the following random variables:%
\begin{equation*}
\begin{split}
  \EMn{X}{n}_0 &= x, \\
  \EMn{X}{n}_{(k+1)/n} &=\EMn{X}{n}_{k/n} +\frac{1}{n}\tilde{V}_0
  \left(\EMn{X}{n}_{k/n}\right)+\frac{1}{\sqrt{n}}\sum_{i=1}^d
  V_i\left(\EMn{X}{n}_{k/n}\right)Z^i_{k+1}.
\end{split}
\end{equation*}%
Then, one can show \cite{KloedenPlaten:1999}\cite{TalayTubaro:1990}
that for an arbitrary $C^4$ function $f$
\begin{equation}\label{eq:EM}
\left\Vert E\left[ f\left(\EMn{X}{n}_1\right) \right]
  -E\left[f\left(
Y(1,x)\right) \right] \right\Vert \leq C_{f}\frac{1}{n}.
\end{equation}
Of course, one needs an algorithm to compute $E\left[f\left(
\EMn{X}{n}_1\right)\right].$
If the $Z_k$ are constructed from
Bernoulli random variables,
$E\left[f\left(\EMn{X}{n}_1\right)\right]$
is a discrete sum, but one would need to do $2^{nd}$ additions,
which can be rather
lengthy when $nd$ is large
(one is then forced to do some Monte-Carlo on a discrete measure).
If the $Z_k$ are normal random variables, one then is forced to do use some
Monte Carlo or quasi-Monte Carlo techniques. When $nd$ is big, quasi-Monte
Carlo method become less effective than Monte-Carlo, but if $nd$ is not too
high, quasi-Monte Carlo method can be very efficient.
\par Another method with the same rate of convergence
appeared in~\cite{LyonsVictoir:2002}, and
is called cubature on Wiener space of degree 3. It is defined with the
following recursive formula:%
\begin{equation*}
\begin{split}
  \CB3{X}{n}_0 &= x,\\
  \CB3{X}{n}_{(k+1)/n} &=\exp\left(\frac{1}{n}V_0
  +\frac{1}{\sqrt{n}}\sum_{i=1}^dZ^i_{k+1}V_i\right)
  \CB3{X}{n}_{k/n}
\end{split}
\end{equation*}%
Such algorithm can be seen as a practical application of the Wong-Zakai
theorem~\cite{ikeda-watanabe}\mbox{\cite{WongZakai:1965}},
when the $Z_k$ are normal random variables.
\\
If $B^{n}_t=(B^{n,1}_t,\dots B^{n,d}_t)\quad(n\in\mathbb{N})$ is the
piecewise linear approximation of the Brownian motion
defined by
\begin{equation*}
  B_{t}^{n}=\left( \left\lfloor nt\right\rfloor +1-nt\right)
  B_{\left\lfloor nt\right\rfloor/n}
  +\left( nt-\left\lfloor nt\right\rfloor\right)
  B_{(\left\lfloor nt\right\rfloor +1)/n},
\end{equation*}
and $Y^{n}$ denotes the solution of the ordinary differential equation
\begin{equation*}
Y_{t}^{n}=x+\int_{0}^{t}V_{0}\left( Y_{s}^{n}\right)
ds+\sum_{i=1}^{d}\int_{0}^{t}V_{i}\left( Y_{s}^{n}\right) dB_{s}^{n,i},
\end{equation*}
then the Wong-Zakai theorem states that $Y^{n}$ converges almost surely to $%
Y^{x}$. It is easy to see that $\CB3{X}{n}_1$ and $Y^{n}_1$ are equal in law,
proving the convergence of the weak algorithm cubature on Wiener space of
degree $3$
(but this argument does not provide the rate of convergence).


\begin{remark}\label{stability}
In the algorithm cubature on Wiener space of degree $3$,
one has to solve numerically ODEs (unless one is lucky and one has a close
form solution!). One possibility is to take its Taylor approximation of
order $1$ for the approximation of $\exp \left( V\right)x$
and we fall back on an Euler scheme. Taking a better approximation (Taylor
approximation of order $2)$ will give a scheme sometimes described as the
Milstein scheme. Not spending enough care on the approximating method of the
ODEs to be solved can result in some catastrophic situations. A general case
where that happens is when the diffusion is almost surely on a subset of ${%
\mathbb{R}}^{n},$ that is, does not fill the whole space. If one has an
approximation scheme which at some time provides an answer outside this set
(which is what happen if one approximates badly the ODEs), the algorithm may
go very wrong or even bug. Increasing $n$ (which is costly) or artificial
techniques can be implemented to solve this problem, while this can be
overcome by taking an appropriately good approximation of the ODEs which
have to be solved (we usually recommend a high order Runge-Kutta scheme, or
an adaptive step size scheme, but this may depend on the particular SDE to
approximate). We will give an example of this problem in
Section~\ref{Sec:num-example}.
\end{remark}
\begin{remark}\label{rem:cub3}
  Random variables which satisfy~\eqref{eq:cub3cond} are easy to find.
  One can take, for a fixed $i$, $Z_{i}^{j}$ to be $d$ independent Bernoulli
  or Gaussian random variables.
  A more elaborate choice of such random variables appeared
  in~\cite{LyonsVictoir:2002}\cite{stroud:1971}.
\end{remark}
\begin{remark}\label{rem:partition}
  Here, we have used the subdivision
  $\left(k/n\right)_{k\in\left\{0,\cdots ,n\right\} }$
  of $\left[ 0,1\right].$
  It is not clear whether taking equal time steps is optimal or not.
  Recently, Kusuoka~\cite{kusuoka:2005:presentation} proved that
  the partitioning into equal time steps
  is optimal when we use the algorithm which we will propose in this paper.
  We do not want to address this problem in this paper, and
  we will always take subdivisions with equal time steps.
\end{remark}

\subsubsection{Higher order}\label{ss:ho}
A way to obtain approximations of higher order is based on the
understanding of more terms in the stochastic Taylor formula~(see
\cite{Castell:1993} and \cite{KloedenPlaten:1999} for example).
When the vector fields $V_{i}$ commute, it is relatively easy to find a
scheme of high order, see \cite{KloedenPlaten:1999} and the
references within. In the general case, one needs to understand how to
approximate weakly the increments of the Brownian motion together with
its first few iterated integrals. This was first successfully done, to
our knowledge,
in~\cite{kusuoka:2001aprx}\cite{LiuLi:2000}\cite{Talay:1990}\cite{Talay:1995}\cite{KusuokaNinomiya:2004}
and then generalized with the method cubature on Wiener
space~\cite{LyonsVictoir:2002}.

\subsection{Romberg Extrapolation\label{ss:Romberg}}

Consider a nice scheme of order $p$, that is, a scheme $\ODp{X}{n}_{k/n}$
such that for smooth $f$, there exists a constant $K_{f}$ such that 
\begin{equation*}
  \left\vert E\left[ f\left( \ODp{X}{n}_{1}\right) \right] -E\left[ f\left(
    Y(1,x)\right) \right] -K_{f}\frac{1}{n^{p}}\right\vert \leq C_{f}\frac{1}{%
    n^{p+1}}.
\end{equation*}%
Then, 
\begin{equation}\label{eq:romberg}
  \frac{2^{p}}{2^{p}-1}E\left[f\left(\ODp{X}{2n}_{1}\right)\right]-\frac{1
  }{2^{p}-1}E\left[f\left(\ODp{X}{n}_{1}\right)\right]
\end{equation}
provides a scheme of order $p+1$. We refer once again
to~\cite{TalayTubaro:1990} for more details and the proof
that the Euler-Maruyama scheme and
its successive Romberg extrapolations are ``nice'' schemes.
Recently, it was proved that our new algorithm presented below is a ``nice''
scheme~\cite{kusuoka:2005:presentation}.

\subsection{A remark on the Monte Carlo method}
\label{ss:RMC}
Let $W$ be a random variable.
When we compute $E[W]$ by Monte-Carlo method with $M$ samples,
we consider a random variable $\left(\sum_{k=1}^MW_k\right)/M$
where $W_i$'s are independent random variables
whose distributions are identical to $W$'s.
We denote this random variable by $\MC(W,M)$.
By virtue of the central limit theorem,
we can consider that $\MC(W,M)$ behaves as a normal random
variable of mean $E[W]$ and variance $\VAR[W]/M$.
\par
Let $\ODp{X}{n}_{1}$ denotes a scheme of order $p$ of the type above.
To calculate $\ODp{X}{n}_{1}$ numerically,
one need to approximate an integral over a $nC(d)$
dimensional space ($C(d)$ denoting a function depending on $d$; for Euler or
Cub3, $C(d)=d$.
As we will see later, $C(d)=d+1$ for our new algorithm).
If one uses the Monte-Carlo method
to approximate this integral, and uses $M$ samples,
the random variable
$\MC\left(f\left(\ODp{X}{n}_{1}\right),M\right)$
is considered.
The situation is summarized as following relations:
\begin{gather}
  \label{eq:DescErr}
  E\left[f(Y(1,x))\right] = E\left[f\left(\ODp{X}{n}_{1}\right)\right]
  + O\left(n^{-p}\right), \\
  \label{eq:IntErr}
  \MC\left(f\left(\ODp{X}{n}_{1}\right), M\right)
  \sim N\left(E\left[f\left(\ODp{X}{n}_{1}\right)\right],
  \frac{\VAR\left[f\left(\ODp{X}{n}_{1}\right)\right]}{M}\right).
\end{gather}
Two types of approximation errors are involved in this calculation.
One is the difference between $E\left[f(Y(1,x))\right]$ and
$E\left[f\left(\ODp{X}{n}_{1}\right)\right]$ and the other
is the difference
between $\MC\left(f\left(\ODp{X}{n}_{1}\right), M\right)(\omega)$
and $E\left[f\left(\ODp{X}{n}_{1}\right)\right]$.
In this paper, we call the former error discretization error and
the latter error integration error.
\eqref{eq:IntErr} shows that we can consider the integration error of
Monte Carlo method to be
a normal random variable of mean $0$ and variance
$\VAR\left[f\left(\ODp{X}{n}_{1}\right)\right]/M$.
\par
Because the difference
between $\VAR\left[f\left(\ODp{X}{n}_{1}\right)\right]$ and
$\VAR\left[f\left(Y(1,x)\right)\right]$ is very small,
we get the following remark.
\begin{remark}\label{rem:MC}
  As long as we use the Monte Carlo method for numerical approximation
  of $E[f(Y(1,x))]$,
  the number of sample points needed to attain the given
  accuracy is independent of the dimension of integration,
  namely the number $n$ of partitions
  and the order $p$ of the approximation scheme.
\end{remark}

\subsection{A remark on the quasi-Monte Carlo method}
Although there are some results which justify the quasi-Monte Carlo method
and give theoretical error with respect to the number $M$ of sample points
and the dimension of the integral domain, those results help little for
error estimation in practice when we apply the quasi-Monte Carlo method to
weak approximation of SDEs
(see~\cite{NinomiyaTezuka:1996} or~\cite{paskov:1997}).
The following observation seems to be widely accepted:
\begin{remark}\label{rem:QMC}
  In contrast to the Monte Carlo case,
  the number of sample points needed by the quasi-Monte Carlo method
  for numerical approximation of $E[f(Y(1,x))]$
  depends heavily on the dimension of integration.
  Smaller the dimension, smaller number of samples are needed.
\end{remark}


The integral that we have to approximate
to obtain $X_{1}^{(\mathrm{ord}\,p),n}$ is on a space of dimension $nC(d)$.
If the numerical method is of high order and $nC\left( d\right) $ is not too
big, one can then use quasi-Monte Carlo with this numerical method to obtain
a very fast algorithm.

Therefore, it seems optimal to look for a (simple) scheme of order
greater than that of the Euler-Maruyama scheme (one), with
$C\left(d\right)$
remaining comparable to
$d$ (i.e. the $C\left( d\right)$
of the Euler-Maruyama scheme).
This is the object of this paper,
where we suggest a new numerical scheme of order $2$, with
$C\left(d\right) =d+1.$
We will show its efficiency by numerically pricing an
Asian option under the Heston model.

\section{Presentation of the new Algorithm}
\label{s:PNA}
We present our new algorithm, of order 2.

\begin{theorem}\label{th:main}
Let $\left( \Lambda _{i},Z_{i}\right) _{i\in \left\{ 1,\cdots ,n\right\} }$
be $n$ independent random variables,
where each $\Lambda _{i}$ is a Bernoulli random
variable independent of $Z_{i}$,
which is a standard $d$-dimensional normal random
variable.
Define $\{\NCUB{X}{n}_{k/n}\}_{k=0,\dots,n}$ to be
a family of random variables 
as follows:
\begin{equation}\label{eq:NewAlgorithm}\begin{split}
    &\NCUB{X}{n}_0 = x,\\
    &\NCUB{X}{n}_{(k+1)/n}=\\
    &\quad\begin{cases}
    \exp\left(\dfrac{V_0}{2n}\right)\exp\left(\dfrac{Z_k^1V_1}{\sqrt{n}}\right)
    \cdots\exp\left(\dfrac{Z_k^dV_d}{\sqrt{n}}\right)
    \exp\left(\dfrac{V_0}{2n}\right)
    \NCUB{X}{n}_{k/n} & \text{if $\Lambda_k = +1,$} \\
    \exp\left(\dfrac{V_0}{2n}\right)\exp\left(\dfrac{Z_k^dV_d}{\sqrt{n}}\right)
    \cdots\exp\left(\dfrac{Z_k^1V_1}{\sqrt{n}}\right)
    \exp\left(\dfrac{V_0}{2n}\right)
    \NCUB{X}{n}_{k/n} & \text{if $\Lambda_k = -1.$}
    \end{cases}
\end{split}\end{equation}
Then, for all $f\in C^\infty_b(\mathbb{R}^N)$,
\begin{equation*}
  \left\vert E\left[f\left(\NCUB{X}{n}_1\right) \right]
  -E\left[ f\left(
  Y(1,x)\right) \right] \right\vert \leq \frac{C_f}{n^2},
\end{equation*}
that is, our new algorithm is of order $2$.
\end{theorem}

A few remarks before all: To compute
\begin{equation*}
  \exp\left(\dfrac{V_0}{2n}\right)\exp\left(\dfrac{Z_k^1V_1}{\sqrt{n}}\right)
  \cdots\exp\left(\dfrac{Z_k^dV_d}{\sqrt{n}}\right)
  \exp\left(\dfrac{V_0}{2n}\right)
  \NCUB{X}{n}_{k/n},
\end{equation*}
one needs to solve $d+2$ ordinary differential
equations.
First along the vector field
$V_0$ from $t=0$ to $t=1/(2n)$
with starting point $\NCUB{X}{n}_{k/n}$, then along
$V_d$ from $t=0$ to $t=Z_k^d/\sqrt{n}$
with starting point the solution of the ODE we have just
solved, and we repeat similar operations $d+2$ times.
One would need an algorithm
to solve this ODE numerically (unless one has a close form solution), and
we, once again, strongly suggest
that one pays 
a lot of attention to the quality of
such algorithm.
\par
One of course will have to use an algorithm to approximate
$E\left[f\left(\NCUB{X}{n}_{1}\right)\right]$, but
this is just a (difficult but classical,
common to Euler algorithm for example) problem of integrating a function on
a finite dimensional space. The simplest but quite effective method is to do
some basic Monte-Carlo simulation of the random variables
$\left( \Lambda_{i},Z_{i}\right) _{i\in \left\{ 1,\cdots ,n\right\} }$.
One could also
simulate the random variables
$\left(\Lambda_i, Z_{i}\right) _{i\in \left\{ 1,\cdots,n\right\} }$
with some quasi-Monte Carlo techniques, or replace the random
variables $Z_{i}$ with some discrete random variables with the right moment up
to order $5$. As this is a very classical problem and common to all the
other probabilistic solutions to our numerical problem, we do not provide
anymore precisions here.
\begin{proof}
The proof is quite classical, so we will not go into details.
The reader should be convinced that the
algorithm is of order $2$ once we show that for $f$ smooth enough,
\begin{equation*}
  \left\vert E\left[f\left(\NCUB{X}{n}_{1/n}\right)\right]
  -E\left[f\left(Y(1/n,x)\right)\right]\right\vert
  \leq \frac{C_f}{n^3}.
\end{equation*}%
The error over $n$ steps, from the Markov property of $Y$,
would then be $n$ times $n^{-3}$.
We consider a smooth function $f$.
First observe that, from the Feynman-Kac theorem,
\begin{equation*}
  \left\vert E\left[ f\left(Y(1/n,x)\right) \right]
  -\left( f(x)+\frac{1}{n}%
  Lf(x)+\frac{1}{2n^{2}}L^{2}f\left( x\right) \right) \right\vert \leq
  C^\prime_{f}n^{-3}.
\end{equation*}%
Developing $L^{2},$ that means
\begin{multline*}
  f(x)+\frac{1}{n}Lf(x)+\frac{1}{2n^{2}}L^{2}f\left( x\right) =
  f(x)+\frac{1}{n}\left( V_{0}+\frac{1}{2}\sum_{i=1}^{d}V_{i}^{2}\right)f(x)\\
  +\frac{1}{2n^{2}}\left( V_{0}^{2}+\frac{1}{2}V_{0}
  \sum_{i=1}^{d}V_{i}^{2}+
  \frac{1}{2}\sum_{i=1}^{d}V_{i}^{2}V_{0}+\frac{1}{4}%
  \sum_{i,j=1}^{d}V_{i}^{2}V_{j}^{2}\right)f(x).
\end{multline*}

Now we need to approximate
$E\left[ f\left(\NCUB{X}{n}_{1/n}\right) \right].$
Using Taylor approximation of the ODEs involved, we quickly see that
the absolute value of
\begin{equation*}
  E\left[
  f\left( \exp \left( \frac{1}{2n}V_{0}\right) \exp \left( \frac{1}{\sqrt{n}}
  Z_{k}^{1}V_{1}\right) \cdots \exp \left( \frac{1}{\sqrt{n}}%
  Z_{k}^{d}V_{d}\right) \exp \left( \frac{1}{2n}V_{0}\right)
  x \right)
  \right]
\end{equation*}%
minus 
\begin{equation*}\begin{split}
  &f(x)+\frac{1}{n}\left( V_{0}+\frac{1}{2}\sum_{i=1}^{d}V_{i}^{2}\right)f(x)\\
  &+\frac{1}{2n^{2}}\left( V_{0}^{2}+\frac{1}{2}V_{0}
  \sum_{i=1}^{d}V_{i}^{2}+
  \frac{1}{2}\sum_{i=1}^{d}V_{i}^{2}V_{0}+\frac{1}{4}\sum_{i=1}^dV_i^4+
  \frac{1}{2}\sum_{i<j}^{d}V_{i}^{2}V_{j}^{2}\right)f(x)
\end{split}\end{equation*}
is bounded by $C^{\prime\prime}_{f}n^{-3}.$
Inverting the order in which the vector
fields are integrated, we obtain that the absolute value of
\begin{equation*}\begin{split}
    E\left[
      f\left(
      \exp \left( \frac{1}{2n}V_{0}\right) \exp \left( \frac{1}{\sqrt{n}}
      Z_{k}^{d}V_{d}\right) \cdots \exp \left( \frac{1}{\sqrt{n}}
      Z_{k}^{1}V_{1}\right) \exp \left( \frac{1}{2n}V_{0}\right)
      x \right)
      \right]
\end{split}\end{equation*}
minus
\begin{equation*}\begin{split}
    &f(x)+\frac{1}{n}\left( V_{0}+\frac{1}{2}\sum_{i=1}^{d}V_{i}^{2}\right)
    f\left( x\right) \\
    &+\frac{1}{2n^{2}}
    \left( V_{0}^{2}+\frac{1}{2}V_{0}\sum_{i=1}^{d}V_{i}^{2}+
    \frac{1}{2}\sum_{i=1}^{d}V_{i}^{2}V_{0}
    +\frac{1}{4}\sum_{i=1}^{d}V_{i}^{4}+
    \frac{1}{2}\sum_{i>j}V_{i}^{2}V_{j}^{2}\right)f(x)
\end{split}\end{equation*}
is bounded by $C^{\prime\prime}_{f}n^{-3}.$ Adding up and dividing by $2$, we
obtain that
\begin{equation*}
  \left\vert E\left[f\left(\NCUB{X}{n}_{1/n}\right)\right]
  -E\left[f\left(Y(1/n,x)\right)\right]\right\vert
  \leq \frac{C^\prime_f + C^{\prime\prime}_f}{n^3}.
\end{equation*}%

\end{proof}

\begin{remark}
Using the results in~ \cite{kusuoka:2001aprx}
and~\cite{kusuoka:2004revisited},
one can show the convergence of the algorithm
with $f$ Lipschitz continuous, under a condition on the vector fields
weaker than H\"{o}rmander condition.
We do not do it here to avoid writing a very technical paper.
\end{remark}

This algorithm could be seen in a non-trivial way as a particular case of
the algorithm cubature on Wiener space of degree $5$. One should also notice
some common features with splitting methods.
\par


\section{Numerical Example: Application to Finance}
\label{Sec:num-example}
In this section, we numerically compare our new algorithm
to the Euler-Maruyama scheme and their Romberg extrapolation.
We calculate the price of an Asian call option with maturity $T$
and strike $K$ written on an asset whose price process $Y_1$
satisfies the
following two factor stochastic volatility model
(Heston model~\cite{heston:1993}):
\begin{equation}\label{eq:heston}
  \begin{split}
    Y_1(t, x) &= x_1+\int_0^t\mu Y_1(s, x)\,ds
    +\int_0^t Y_1(s, x)\sqrt{Y_2(s,x)}\,dB^1(s), \\
    Y_2(t, x) &= x_2+\int_0^t\alpha
    \left(\theta - Y_2(s, x)\right)\,ds +\int_0^t\beta
    \sqrt{Y_2(s, x)}\,dB^2(s),
  \end{split}
\end{equation}
where $x = (x_1, x_2) \in(\mathbb{R}_{>0})^2$,
$(B^1(t), B^2(t))$ is a $2$-dimensional standard Brownian motion, and
$\alpha$, $\theta$, $\mu$ are some positive coefficients such that
$2\alpha\theta-\beta^2>0$ to ensure the existence and uniqueness of a
solution to our SDE \cite{feller:1950}.
The payoff of this option is $\max\left(Y_3(T,x)/T-K, 0\right)$,
where
\begin{equation}\label{eq:heston2}
  Y_3(t,x) = \int_0^tY_1(s,x)\,ds.
\end{equation}
The price of this option becomes
$D\times E\left[\max\left(Y_3(T,x)/T - K,\, 0\right)\right]$ where $D$ is
the appropriate discount factor.
We set $T=1$, $K=1.05$, $\mu = 0.05$,
$\alpha = 2.0$, $\beta = 0.1$, $\theta = 0.09$, and
$(x_1, x_2) = (1.0, 0.09)$.
We ignore $D$
in this experiment.
Let $Y(t,x) = \big{.}^t\!\left(Y_1(t,x), Y_2(t,x), Y_3(t,x)\right)$.
We transform the SDEs~\eqref{eq:heston} and~\eqref{eq:heston2}
into a Stratonovich form SDE:
\begin{equation}\label{eq:heston3}
    Y(t,x) = \sum_{i=0}^2\int_0^tV_i(Y(s,x))\circ dB^i(s),
\end{equation}
where
\begin{equation}\begin{split}
    V_0\left({}^t\!\left(y_1, y_2, y_3\right)\right) &=
    \bigg{.}^t\!
    \left(y_1\left(\mu - \frac{y_2}{2}\right),\,
    \alpha(\theta - y_2)-\frac{\beta^2}{4},\,
    y_1\right) \\
    V_1\left({}^t\!\left(y_1, y_2, y_3\right)\right) &=
    \Big{.}^t\!
    \left(y_1\sqrt{y_2},\, 0,\, 0\right) \\
    V_2\left({}^t\!\left(y_1, y_2, y_3\right)\right) &=
    \Big{.}^t\!
    \left(0,\, \beta\sqrt{y_2},\, 0\right).
\end{split}\end{equation}

\subsection{Implementation of the algorithm}
We apply the algorithm which we introduced in Section~\ref{s:PNA} to
this problem.
\subsubsection{Solutions of the ODEs}
We can easily get $\exp\left(sV_1\right)$ and $\exp\left(sV_2\right)$
$(s\in\mathbb{R})$ as follows:
\begin{equation}\label{eq:expV1}
  \begin{split}
    \exp\left(sV_1\right){}^t\!\left(y_1,y_2,y_3\right) &=
    \Big{.}^t\!
    \left(y_1e^{s\sqrt{y_2}}, y_2, y_3\right), \\
    \exp\left(sV_2\right){}^t\!\left(y_1, y_2, y_3\right) &=
    \Bigg{.}^t\!
    \left(y_1, \left(\frac{\beta s}{2}+\sqrt{y_2}\right)^2, y_3\right).
  \end{split}
\end{equation}
As there exists no closed form solution to  $\exp\left(sV_0\right)$,
we are forced to use an approximation and we choose:
\begin{equation}\label{eq:expV0}
    \exp\left(sV_0\right){}^t\!\left(y_1, y_2, y_3\right)
    = \big{.}^t\!\left(g_1(s),g_2(s), g_3(s)\right),
\end{equation}
where
\begin{equation}\label{eq:expV0ii}
  \begin{split}
    g_1(s) &= y_1\exp\left(\left(\mu-\frac{J}{2}\right)s
    +\frac{y_2-J}{2\alpha}\left(e^{-\alpha s}-1\right)\right), \\
    g_2(s) &= J+\left(y_2 - J\right)e^{-\alpha s}, \\
    g_3(s) &= y_3 + \frac{y_1\left(e^{As}-1\right)}{A}
    + O\left(s^3\right),
    \\
      J&=\theta - \frac{\beta^2}{4\alpha},\quad
      \text{and}\quad
      A=\mu-\frac{y_2}{2}.
  \end{split}
\end{equation}
The error compared to the true solution is $O\left( t^{3}\right) $ in small
time $t$,
creating an additional error of $O\left( n^{-3}\right) $ at every step
of the algorithm, but as the error of our scheme at every step was also $%
O\left( n^{-3}\right) ,$ taking the above approximation of $\exp \left(
sV_{0}\right) $ does not alter the convergence rate of the algorithm.
\par
Following the same discussion, it is easy to see that we have to
approximate $\exp(tV_0)$ in such a way that the order of the produced error
is $O(t^4)$ when we use Romberg extrapolation,
which we introduced in~\ref{ss:Romberg}, together with the algorithm.
In this experiment,
we approximate $g_3(s)$ by the traditional order $4$ Runge-Kutta
method when we use Romberg extrapolation.
\par
Here, we see that
one of the advantages of this algorithm over the Euler-Maruyama
scheme is the one we mentioned in Remark~\ref{stability}.
When we apply the Euler-Maruyama scheme to this
process~\eqref{eq:heston},
it may
happen that the square volatility process $(Y_{2})_{k}^{(\mathrm{EM}),n}$
becomes negative, and the algorithm then fails at the next step (as we will
have to take its square root). On the other hand,
equations~\eqref{eq:expV1} and~\eqref{eq:expV0ii} show
that our new algorithm does not share this problem.
There exists a way of avoiding this problem with the
Euler-Maruyama scheme~\cite{diop:2004}.

\subsubsection{A remark on general implementation}
In general, it is not always possible to obtain
the closed form solution to $\exp(sV_i)$.
Even in such cases,
it is not difficult to implement our new algorithm.
All we have to do is to find an approximation
of $\exp(sV_0)$ whose error is $O(s^3)$ and
approximations of $\exp(sV_i),\ (i\neq 0)$ whose
errors are $O(s^6)$.
This can be achieved by Runge-Kutta like methods
and we can find some examples of them in \cite{butcher:1987}.
\par
We remark that when we use Romberg extrapolation together,
we have to approximate $\exp(sV_0)$ with $O(s^4)$ error
and $\exp(sV_i)\, (i\neq 0)$ with $O(s^7)$ error.

\subsubsection{Application of the quasi-Monte Carlo method}
Our new algorithm has the virtue that
the application of the quasi-Monte Carlo method to this algorithm
is possible in a straight forward way,
once we embed
$(\Lambda_i, Z_i)_{i\in\{1,\dots,n\}}$ into $[0,1)^{n(d+1)}$.
This is an advantage of the algorithm
over algorithms proposed
in~\cite{ninomiya:2001a},
\cite{ninomiya:2003}, and~\cite{KusuokaNinomiya:2004}
which also enable us to proceed higher order weak approximation.


\subsection{Comparison to Euler-Maruyama scheme}
\begin{figure}[h]
  \centering{\includegraphics{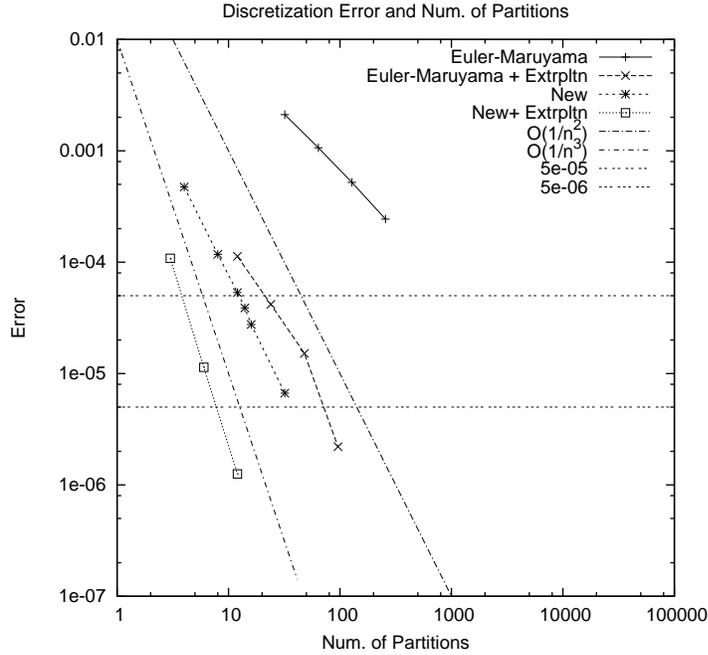}}
 \caption{Error coming from the discretization}
 \label{fig:disc-err}
\end{figure}
We compare numerically our new algorithm
to the Euler-Maruyama scheme with and without Romberg
extrapolation. Such methods involve, as we saw, approximation of an integral
over a finite dimensional space; we will do these approximations using the
Monte Carlo method and the quasi-Monte Carlo method.
\par
There are many studies on acceleration of Monte Carlo
methods~\cite{Glasserman:2004}
but we choose the crude Euler-Maruyama scheme with and without Romberg
extrapolation as only competitors by the following reasons:
\begin{enumerate}
\item Only our new algorithm and the Euler-Maruyama scheme are
  very universal and applicable easily
  to any type of problems described in subsection~\ref{ss:problem}.
\item Almost all of variance reduction techniques which we can apply to
  the Euler-Maruyama scheme are also applicable to our new algorithm.
\end{enumerate}
These are important advantages of our new algorithm.  Many existing
algorithms lack one or both of these properties.  For example,
in~\cite{LapeyreTemam:2001}, they proposed the trapezoidal algorithm
which accelerates Monte Carlo pricing of Asian option price.  But this
algorithm works only for pricing of Asian option written on one
dimensional diffusion.  There are many such type of problem-specific
algorithms and we exclude them, because in this paper we focus on
universal algorithms which work for any weak approximation problem of
any diffusion processes defined by~\eqref{StSDE}.

\par
In this experiment, we consider
$$E\left[\max\left(Y_3(T,x)/T - K,\, 0\right)\right]
=6.0473907415
\times 10^{-2}$$
which is obtained by our new algorithm with extrapolation,
quasi-Monte Carlo,
$n=96 + 48$, and $M=8.0\times 10^9$.
\subsubsection{Discretization Error}

Figure~\ref{fig:disc-err} shows the relation between the number of
partitions in our discretization of the interval $\left[ 0,1\right] $ ($n$
in the description of the algorithm) and the error of the algorithms. We
observe that to achieve $10^{-4}$ accuracy, our new method with Romberg
extrapolation requires $n=6$, our new method needs $n=12$, while the
Euler-Maruyama scheme with Romberg extrapolation needs $n=24$, and the
simple Euler-Maruyama scheme needs $n\geq 2000$.
In all algorithms,
consumed time is proportional to $n\times M$, where $M$ is the number
of sample points.

\subsubsection{Convergence Error from Monte Carlo}
\label{ss:CEMC}
\begin{figure}[h]
  \centering{\includegraphics{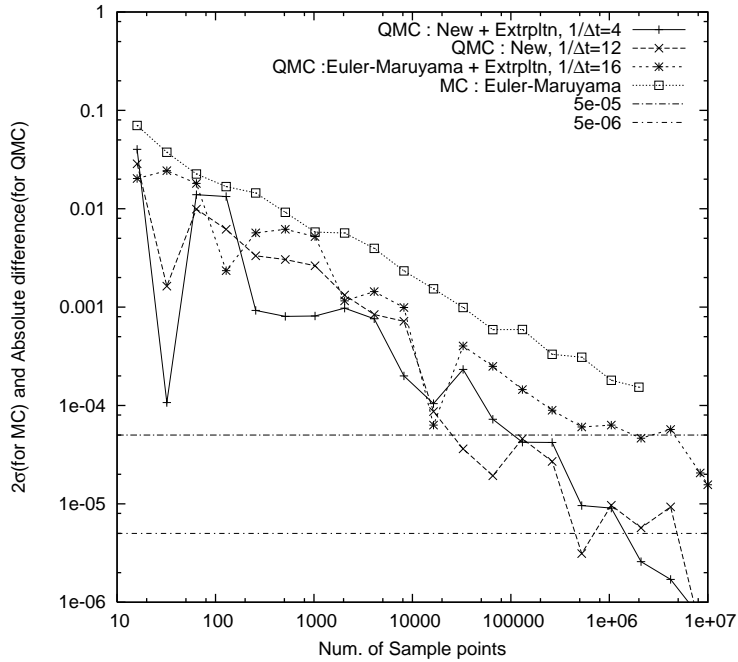}}
  \caption{Convergence Error from quasi-Monte Carlo and Monte Carlo}
  \label{fig:QMC_conv}
\end{figure}
\noindent
We have already mentioned in~\ref{ss:RMC} that the
convergence performance of the Monte Carlo method
is independent of the number of partitions.
We can see in
Figure~\ref{fig:QMC_conv}
that
in this experiment this statement holds.
This figure also shows that to achieve $10^{-4}$ accuracy with
95\% confidence level ($2\sigma$) by using Monte Carlo method,
we need over $10^8$ sample points.
We can also see in this figure
that the Monte Carlo errors which come from algorithms boosted by
the Romberg extrapolation become greater than those of the original
algorithms.
\subsubsection{Convergence Error from quasi-Monte Carlo and Monte Carlo}
Figure~\ref{fig:QMC_conv} also shows that the performance of
the convergence of the quasi-Monte Carlo method depends on the
number $n$ of partitions and on the algorithms.
Figure~\ref{fig:QMC_conv} seems to show that the quasi-Monte Carlo
method outperforms the Monte Carlo method specially when used
with our new algorithm
and that the algorithm needs $2\times 10^5$ sample points
for $10^{-4}$ accuracy,
the algorithm with extrapolation $2\times 10^5$ sample points,
and Euler-Maruyama with extrapolation $5\times 10^6$ sample points when
we use the quasi-Monte Carlo method.
\subsubsection{Performance comparison with respect to consumed time}
\begin{table}[h]
  \begin{center}
    \caption{\#Partition, \#Sample, and CPU time required for
      $10^{-4}$ accuracy.}\label{table:table}
    \begin{tabular}{|l|ccc|}
      \hline
      Method & \#Partition & \#Sample & CPU time (sec) \\
      \hline
      \hline
      E-M + MC & $2000$ & $10^8$ & $1.72\times 10^5$ \\
      \hline
      E-M + Extrpltn + MC & $16 + 8$ & $10^8$ & $2.06\times 10^3$ \\
      \hline
      New + MC & $12$ & $10^8$ & $1.24\times 10^3$ \\
      \hline
      New + Extrpltn + MC & $4+2$ & $10^8$ & $6.2\times 10^2$ \\
      \hline
      E-M + Extrpltn + QMC & $16 + 8$ & $5\times 10^6$ & $1.28\times 10^2$ \\
      \hline
      New + QMC & $12$ & $2\times 10^5$ & $3.3$ \\
      \hline
      New + Extrpltn + QMC & $4+2$ & $2\times 10^5$ & $1.73$ \\
      \hline
    \end{tabular}
  \end{center}
\end{table}
\noindent
The elapsed time of each method required for $10^{-4}$ accuracy
is shown in Table~\ref{table:table}.
We find in this table that our new algorithm with Romberg extrapolation
and the quasi-Monte Carlo
method provides the fastest calculation.
Our new algorithm with Romberg extrapolation and quasi-Monte Carlo
is about $80$ times faster than Euler-Maruyama scheme
with Romberg extrapolation and quasi-Monte Carlo.
We also see that even without Romberg extrapolation,
our new algorithm is still faster than any
boosted Euler-Maruyama method.
\par
At last we would like to mention
Remark~\ref{rem:MC} and Remark~\ref{rem:QMC} again.
The remarkable performance of our new algorithm is
closely related to the property of the quasi-Monte Carlo method
noted in Remark~\ref{rem:QMC}.


\bibliographystyle{amsplain}
\bibliography{ninomiya}
\end{document}